\documentclass[11pt]{amsart}
\usepackage{amsfonts}
\usepackage{amsmath}
\usepackage{amsthm}
\usepackage{url}
\usepackage[all]{xy}
\usepackage{latexsym}
\usepackage{amssymb}
\usepackage[cp850]{inputenc}
\usepackage{amsfonts}
\usepackage{graphicx}

\usepackage[abs]{overpic}

\usepackage[usenames]{color}

\newtheorem{example}{Example}

\newtheorem*{claim}{Claim}

\newtheorem*{defi*}{Definition}			
\newtheorem*{bei*}{Example}
\newtheorem*{sat*}{Theorem}				
\newtheorem*{kor*}{Corollary}
\newtheorem*{rmk*}{Remark}				
	
\newtheorem*{quest*}{Question}

\let\ssection=\section
\renewcommand{\section}{\setcounter{equation}{0}\ssection}

\newtheorem*{namedtheorem}{\theoremname}
\newcommand{\theoremname}{testing}

\theoremstyle{remark}

\newtheorem*{namedtheoremr}{\theoremnamer}
\newcommand{\theoremnamer}{testing}

\newcommand{\BC}{\mathbb C}			\newcommand{\BH}{\mathbb H}
\newcommand{\BR}{\mathbb R}			
\newcommand{\BN}{\mathbb N}			\newcommand{\BQ}{\mathbb Q}
\newcommand{\BS}{\mathbb S}

\newcommand{\CE}{\mathcal E}

		\newcommand{\CT}{\mathcal T}
		\newcommand{\CV}{\mathcal V}

\newcommand{\D}{\partial}

\DeclareMathOperator{\PSL}{PSL}		
\DeclareMathOperator{\Id}{Id}		
\DeclareMathOperator{\Isom}{Isom}	

\DeclareMathOperator{\Ker}{Ker}

\newcommand{\comment}[1]{}

\begin{document}

\title[]{Discrete groups without finite quotients}
\author{Tommaso Cremaschi}
\address{Department of Mathematics, Boston College}
\email{cremasch@bc.edu}
\author{Juan Souto}
\address{IRMAR, Universit\'e de Rennes 1}
\email{juan.souto@univ-rennes1.fr}
\date{\today}

\thanks{T.C. was partially supported by PICS \#7734.}

\begin{abstract}
We construct an infinite discrete subgroup of the isometry group of $\BH^3$ with no finite quotients other than the trivial group.
\end{abstract}

\maketitle

It is well-known that every finitely generated linear group is residually finite \cite{Selberg}. Finite generation is definitively necessary, as it is already made apparent by the group $\BQ$. However, people working with Kleinian groups, that is with discrete groups of isometries of hyperbolic spaces $\BH^n$, might find examples as $\BQ$ to be kind of pathological. In fact, it is well-known that discreteness of a group of isometries of hyperbolic space imposes non-trivial algebraic conditions. For example, centralisers of infinite order elements in discrete Kleinian groups are virtually abelian. Or, more to the point, while $\PSL_2\BQ\subset\PSL_2\BR\subset\Isom(\BH^2)$ is simple \cite{Dickson,Moore}, it is easy to see, using small cancellation arguments, that there are no infinite, simple, and discrete subgroups of $\Isom(\BH^n)$ (compare with \cite{Delzant,Delzant-Gromov}). Also, Kleinian groups are mostly studied in low dimensions, and in that setting further algebraic restrictions do arise. For instance, the classification of 2-dimensional orbifolds implies that all discrete subgroups of $\Isom\BH^2$ are residually finite. 

The goal of this note is to present examples of discrete subgroups of $\Isom\BH^3$ which fail to be residually finite. In fact, they don't have any non-trivial finite quotients whatsoever. 

\begin{example}\label{ex3}
There is an infinite discrete subgroup $G\subset\Isom(\BH^3)$ without finite non-trivial quotients. 
\end{example}

As we just said, having no finite quotients, the group $G$ in Example \ref{ex3} clearly fails to be residually finite. Examples of discrete, non-residually finite subgroups of $\Isom(\BH^3)$ have been previously constructed by Agol \cite{Agol}. Both Agol's examples and the group in Example \ref{ex3} have torsion. We present next an example, a variation of Agol's example, showing that there are also torsion free discrete non-residually finite subgroups of $\Isom(\BH^3)$:

\begin{example}\label{ex1}
There is a torsion free discrete subgroup $G\subset\Isom(\BH^3)$ which is not residually finite.
\end{example}

The remaining of this note is devoted to discuss these two examples.

\section*{Discussion of the examples}

In the course of our discussion we feel free to use standard facts of hyperbolic geometry as one might find in classical texts such as \cite{Maskit,Ratcliffe}. It will also be convenient to see our groups as fundamental groups of infinite, locally finite, graph of groups. We refer to standard texts like \cite{Serre} for basic facts about graphs of groups.

\subsection*{Example 1.}
We give an algebraic description of a group $G$, then we prove that it has no finite quotients, and we finally show that it is isomorphic to a discrete subgroup of $\PSL_2\BC$.

Let $T$ be the maximal rooted binary tree. Denote by $\CV$ and $\CE$ the sets of vertices and edges respectively, let $*$ be the root of $T$ and, for $v\in\CV$, let $\vert v\vert\in\BN$ be the distance from $v$ to $*$. We orient the edges of $T$ so that they point to the root and for $e\in\CE$ we let $e^+$ be its terminal vertex. Given a vertex $v\in\CV$ with $\vert v\vert\ge 1$ let $e_0(v)$ the edge leaving $v$ and pointing out of $v$ and label the two edges pointing into $v$ by $e_1(v)$ and $e_2(v)$. 

Consider from now the group 
$$G=\left\langle\{g_e\vert e\in\CE\}\middle\vert\left\{g_{e_0(v)}^{3+\vert v\vert},\ g_{e_0(v)}g_{e_1(v)}^{-1}g_{e_2(v)}^{-1}\middle\vert v\in\CV\text{ with }\vert v\vert\ge 1\right\}\right\rangle.$$
The group $G$ also admits a description as the fundamental group 
$$G=\pi_1(\CT)$$
of a graph of groups $\CT$ with underlaying graph $T$ 
with vertex groups
$$G_v=\left\langle g_{e_0(v)},g_{e_1(v)},g_{e_2(v)}\middle\vert g_{e_0(v)}^{3+\vert v\vert},g_{e_1(v)}^{4+\vert v\vert},g_{e_2(v)}^{4+\vert v\vert},g_{e_0(v)}g_{e_1(v)}^{-1}g_{e_2(v)}^{-1}\right\rangle$$
if $v\neq *$, with
$$G_*=\left\langle g_{e_1(*)},g_{e_2(*)}\middle\vert g_{e_1(*)}^{4},g_{e_2(*)}^{4}\right\rangle,$$
and with edge groups 
$$G_e=\left\langle g_e\middle\vert g_e^{4+\vert e^+\vert}\right\rangle.$$

We are going to think of the group $G$ as the nested union of a sequence of subgroups. The easiest way to describe these subgroups is as the fundamental groups
$$G^n=\pi_1(\CT^n)$$
of the subgraph of groups $\CT^n\subset \CT$ corresponding to the ball of radius $n-1$ around the root $*$. Alternatively, $G^n$ is the subgroup of $G$ generated by all those elements $g_{e_0(v)}$ with $\vert v\vert\le n$. We have
$$G^1\subset G^2\subset G^3\subset\ldots,\ \ G=\bigcup_{n\in\mathbb N}G^n.$$

\begin{claim}
The group $G^n$ is generated by the set $S^n=\{g_{e_0(v)}\text{ with }\vert v\vert=n\}$.
\end{claim}
\begin{proof}
Since $G^n$ is generated by $S_1\cup S_2\cup\dots\cup S_n$ and hence by $G^{n-1}\cup S^n$, we can argue by induction on $n$. Therefore, it suffices to prove that $S^{n-1}$ is contained in the group generated by $S^n$. Well, given $v\in\CV$ with $\vert v\vert =n-1$ let $v_1$ and $v_2\in\CV$ be the initial vertices of the edges $e_1(v)$ and $e_2(v)$. Given that $e_1(v)=e_0(v_1)$ and $e_2(v)=e_0(v_2)$ the presentation of the group gives us:
$$g_{e_0(v)}=g_{e_0(v_2)}g_{e_0(v_1)}$$
The claim then immediately follows.
\end{proof}

We are now ready to prove that $G$ has no finite quotients:

\begin{claim}
$G$ has no finite non-trivial quotients.
\end{claim}
\begin{proof}
Let $H$ be a finite group and $\pi:G\to H$ be a homomorphism. We need to prove that $\pi$ is trivial. Now, denote by $\vert H\vert$ the order of $H$, let $k\in\BN$ be arbitrary, and let $n>k$ with 
$$3+n\equiv 1\mod(\vert H\vert).$$
From the presentation of the group $G$ we get that the elements in the set $S^{n}=\{g_{e_0(v)}\text{ with }\vert v\vert=n\}$ have order $3+n$ in $G$. Since $\vert H\vert$ and $3+n$ are prime to each other, we get that the elements in $S^n$ are killed by $\pi$:
$$S^{n}\subset\Ker(\pi).$$
Now, from the previous claim we also get that the subgroup $G^{n}$ is contained in the kernel of $\pi$, meaning that also $G^k\subset\ker(\pi)$. Since $G$ is the union of the subgroups $G^k$ we get that $G\subset\ker(\pi)$. We have proved the claim.
\end{proof}

All is left now is to show that $G$ admits discrete and faithful representations into $\PSL_2\BC$. We are going to construct the desired representation as a limit of representations of the groups $G^n$. 

\begin{claim}
For each $n\ge 0$ there is a discrete and faithful representation
$$\rho_n:G^n\to\PSL_2\BC$$
such that the restriction of $\rho_n$ to $G^{n-1}$ agrees with $\rho_{n-1}$ for every $n\ge 1$.

Moreover, the limit set of the group $\rho_n(G_v)$ bounds a round disk in the discontinuity domain of $\rho_n(G^n)$ for every vertex group $v$ with $\vert v\vert=n$.
\end{claim}

The final assertion of the Claim serves to be able to argue by induction -- the real point is the first assertion because it allows us to define
$$\rho:G\to\PSL_2\BC$$
satisfying $\rho(g)=\rho_n(g)$ if $g\in G^n$. This representation is faithful because each one of the representations $\rho_n$ is. The same argument holds true for discreteness because already the first group $\rho(G_1)$ is non-elementary \cite{Benedetti-Petronio}.

All that it is left is to prove the claim.

\begin{proof}[Proof of the claim]
We will argue by induction. First note that 
$$G^0=G_*=\left\langle g_{e_1(*)},g_{e_2(*)}\middle\vert g_{e_1(*)}^{4},g_{e_2(*)}^{4}\right\rangle$$
is isomorphic to the $(4,4,\infty)$-triangle group:
$$H_0=\langle a,b\vert a^{4}, b^{4}\rangle.$$
We can thus take $\rho_0:G_0\to\PSL_2\BC$ to be the standard fuchsian representation and take the said disk to be any one of the two connected components of the discontinuity domain of $\rho_0(G_0)$. 

Suppose that the claim holds true for $n-1$. For each vertex $v\in\CV$ with $\vert v\vert=n-1$ let $\Delta_v\subset\BS^2=\D_\infty\BH^3$ be the disk in the discontinuity domain bounded by the limit set of $\rho_{n-1}(G_v)$. Note that $\Delta_v$ is $\rho_{n-1}(G_v)$-invariant and that no translate of $\Delta_v$ under $\rho_{n-1}(G^{n-1})$ meets $\Delta_w$ for another vertex $w\neq v$ with $\vert w\vert=n-1$. The disk $\Delta_v$ is the boundary at infinity of a hyperbolic half-space $H_v$ - let $\Delta'_v$ be the hyperbolic plane bounding $H_v$.

Suppose now that we have a vertex $z\in\CV$ with $\vert z\vert=n$ and denote by $z^+$ the terminal vertex of the edge $e=e_0(z)$. We identify the edge group $G_{e_0(z)}$ with the corresponding subgroup of the vertex group $G_{z^+}$. The group $\rho_{n-1}(G_{e_0(z)})$ 
is cyclic and has a unique fixed points $\alpha\in\Delta_v$ and $\alpha'\in\Delta'_v$. For $L>0$ let $x\in[\alpha',\alpha)$ be the point at distance $L$ from $\alpha'$ and let $D$ be the hyperbolic plane containing $x$ and perpendicular to the ray $[x,\alpha)$. The plane $D$ is 
$\rho_{n-1}(G_{e_0(z)})$-invariant. 

Consider now the edge group $G_{e_0(z)}$ as a subgroup of $G_z$. The action of $G_{e_0(z)}$ on $D$ via $\rho_{n-1}$ extends to a discrete action of the triangle group $G_z$ and thus to an action of $G^{(n-1)}*_{G_{e_0(z)}}G_z$. 

Proceeding in this way with all vertices $z\in\CV$ with $\vert z\vert=n$ we then get a representation
$$\rho_L:G^n\to\PSL_2\BC$$
depending on the parameter $L$. By construction all of these constructions extend the representation $\rho_{n-1}$ and it follows from the Klein-Maskit combination theorem that for $L$ large enough the representation $\rho_L$ is discrete and satisfies the additional desired claim.
\end{proof}

This concludes the discussion of Example \ref{ex3}.

\subsection*{Example 2.} Let $T$ be a once holed torus and let $\alpha$ and $\beta$ be simple curves in $T$ which intersect (transversely) in a single point. Let also $U$ be a regular neighborhood of $\beta\times\{0\}$ in the 3-manifold $T\times[-1,1]$, $\mu$ the meridian of $U$, and $\beta'\subset\D U$ the longitude of $U$ isotopic to $\beta\times\{1\}$ in $T\times[-1,1]\setminus U$. Finally, for $n\in\BN$ let $M_n$ be manifold obtained from $T\times[-1,1]\setminus U$ by Dehn-filling the curve $n \mu+ \beta'$. The curve $\beta'$ intersects the new meridian $n$ times, which means that it represents the n-th power of the soul of the new solid torus. We get thus the presentation
$$\pi_1(M_n)=\pi_1(T)*_{\langle\beta\rangle}\pi_1(U)\simeq\langle a,b,c\vert b=c^n\rangle$$
with respect to which the curve $\D T\times\{0\}$ corresponds to the conjugacy class of $[a,b]=[a,c^n]$.

Now, the pair $(M_n,\D T\times[-1,1])$ is a pared manifold, which means that the group $\pi_1(M_n)$ admits a geometrically finite representation 
$$\rho_n:\pi_1(M_n)\to \PSL_2\BC=\Isom_+(\BH^3)$$
with $\rho_n([a,b])=g\in\PSL_2\BC$ where $g(z)=z+1$. 

For $t\in\BR$ let $h_t\in\PSL_2\BC$ be the parabolic element $h_t(z)=z+ti$. Now, a standard combination argument implies that for a sufficiently fast growing sequence $t_n\to\infty$, the group 
$$G=\left\langle \cup_nh_{t_n}\rho_n(\pi_1(M_n))h_{t_n}^{-1}\middle\vert n\in\BN\right\rangle$$
generated by the union of the groups $h_{t_n}\rho_n(\pi_1(M_n))h_{t_n}^{-1}$ is discrete.

Notice now that for each $n$ the element $g\in G$ can be written as
$$g=h_t g h_t^{-1}=h_t\rho_n([a,b])h_t^{-1}=h_t\rho_n([a,c^n])h_t^{-1}.$$
It follows that if $H$ is an arbitrary finite group and if $\pi:G\to H$ is a homomorphism then 
\begin{align*}
\pi(g)&=\pi(h_t)[\pi(\rho_n(a),\rho_{\vert H\vert}(c)^{\vert H\vert}]\pi(h_t)^{-1}\\
&=\pi(h_t)[\pi(\rho_n(a),\Id_H]\pi(h_t)^{-1}\\
&=\Id_H
\end{align*}
where $\vert H\vert$ is the order of $H$. Having proved that $g\in G$ belongs to the kernel of every homomorphism to a finite group, we have showed that $G$ is not residually finite. This concludes the discussion of Example \ref{ex1}.

\end{document}